\documentstyle[11pt]{article}
\newtheorem{thm}{Theorem}[section]
\newtheorem{cor}[thm]{Corollary} 
\newtheorem{theorem}{Theorem}[section]
\newtheorem{definition}[theorem]{Definition}
\def\eps{\epsilon}
\def\epsbr{\overline{\epsilon}}
\def\xbf{{\bf x}\,}
\def\Abf{{\bf A}\,}
\def\Abh{\hat{\bf A}\,}
\def\Bbf{{\bf B}\,}

\def\Dbf{{\bf D}\,}
\def\Ebf{{\bf E}\,}
\def\Fbf{{\bf F}\,}

\def\Ibf{{\bf I}\,}
\def\Kbf{{\bf K}\,}
\def\Lbf{{\bf L}\,}
\def\Mbf{{\bf M}\,}
\def\Pbf{{\bf P}\,}
\def\Qbf{{\bf Q}\,}
\def\Sbf{{\bf S}\,}
\def\Tbf{{\bf T}\,}
\def\Ubf{{\bf U}\,}
\def\Vbf{{\bf V}\,}
\def\Xbf{{\bf X}\,}
\def\ebf{{\bf e}\,}
\def\ubf{{\bf u}\,}

\date{ August 1994}
\title{
Generalized Epsilon-Pseudospectra
\thanks{
This work was motivated by
discussions with Wolfgang Kerner on  magneto-hydrodynamic 
pseudospectra.
The author thanks Satish Reddy for introducing him to 
$\epsilon$-pseudospectra, and 
for providing preprints.
The author also thanks Satish Reddy and Nick Trefethen 
for critically reading this manuscript;    
their suggestions have considerably improved the presentation.
J.~Demmel's comments were also valuable. 
This work was performed under
U.S. Dept. of Energy Grant DE-FG02-86ER53223.}}

\author{Kurt S.~Riedel\thanks{Courant Institute of Mathematical
Sciences, New York University, New York, New York 10012-1185
({\tt riedel@cims.nyu.edu}).} }

\begin{document}

\maketitle

\begin{abstract}
We generalize $\epsilon$-pseudospectra and the associated computational
algorithms to the generalized eigenvalue problem. 
Rank one perturbations are used to determine the $\epsilon$-pseudospectra.

\

\end{abstract}


\thispagestyle{plain}
\markboth{K.S.~RIEDEL}{GENERALIZED $\epsilon$-PSEUDOSPECTRA}

\section{$\epsilon$-pseudospectra of the eigenvalue problem}

Normal matrices have complete sets of {\it orthonormal} eigenvectors,
and therefore the spectral decomposition is useful in  studying
the properties of normal operators. In contrast, the eigenvectors of
non-normal matrices can be nearly linearly dependent, and the 
eigenvalue problem may be highly ill-conditioned. Thus
additional concepts and analysis techniques are useful in examining 
non-normal operators.
$\epsilon$-pseudospectra, introduced by L.N.~Trefethen, have proven to be a
powerful tool in the analysis of non-normal operators [10,11]. 
In examining the
transient behavior of the Orr-Somerfeld operator, Reddy et al. calculated
the $\epsilon$-pseudospectra, the numerical range and the maximum transient
growth [8]. The corresponding analysis for resistive magneto-hydrodynamics
(MHD) generates a
generalized eigenvalue problem [2,5,6,9]. 
In this note, we generalize
$\epsilon$-pseudospectra to the generalized eigenvalue problem.

We begin by reviewing $\epsilon$-pseudospectra for the standard eigenvalue problem
for the $n \times n$ matrix $\bf A$. We use the standard 2-norm,
$||{\bf v}||^2 \equiv {\bf v}^{*}{\bf v}$.
We denote 
the spectrum of $\bf A$ by $\Lambda({\bf A})$ and 
the resolvent set of $\bf A$ by $\rho({\bf A})$.

\begin{definition}
{\rm (Trefethen [10,11]).} {\it
Let $\epsilon \geq 0$ be given. A complex number $z$ is in the
$\epsilon$-pseudospectrum of ${\bf A}$,
which we denote by $\Lambda_{\epsilon}( {\bf A})$,
if one of the following equivalent conditions is satisfied:

(i) the smallest singular value of ${\bf A}-z{\bf I}$ is
less than or equal to $\epsilon$,

(ii) $\exists {\bf u}\in C^n$ such that
$||{\bf u}||^2 = 1$ and 
$ ||({\bf A}-z{\bf I}) {\bf u}||^2 \leq \epsilon^2$,

(iii) $z \in \Lambda ({\bf A})$ or
$z \in \rho ({\bf A})$ and $\exists {\bf u} \in C^n$
such that $||{\bf u}||^2 = 1$
and ${\bf u}^{*} ({\bf A}-z{\bf I})^{-1 *}$
$({\bf A}-z{\bf I})^{-1} {\bf u} \geq 1/
\epsilon^2$,

(iv) $z$ is in the spectrum of
${\bf A}+ \epsilon  {\bf E}$:
$({\bf A} + \epsilon{\bf E}) {\bf u} = z{\bf u} $,
where the matrix ${\bf E}$ satisfies
$\parallel {\bf E} \parallel \leq 1$,

(v) $z$ is in the spectrum of ${\bf A} + \epsilon {\bf v}_2 {\bf v}_1^{*}$,
where
$\parallel {\bf v}_1 \parallel \leq1$ and
$\parallel {\bf v}_2 \parallel \leq1$,

(vi) $z$ is in the spectrum of ${\bf A} - ({\bf Av} - z {\bf v}) {\bf v}^{*}$
where
$\parallel {\bf v} \parallel \leq1$ and
$\parallel {\bf Av} - z {\bf v} \parallel \leq \epsilon$.
}
\end{definition}

The equivalence of (i), (ii), (iii) and (iv) is given in Refs.~[8,10,11].
Conditions (v) and (vi) are new, equivalent definitions of 
$\eps$-pseudospectra.
(vi) $\Rightarrow$ (v) $\Rightarrow$ (iv) is trivial, as is
(vi) $\Longleftrightarrow$ (ii).
The original proof that (ii)
$\Rightarrow$ (iv)  implicitly contained  the proof that (ii) $\Rightarrow$
(vi) since the constructed ${\bf E}$ matrix is given by (vi). \hfill $\Box$

The stated definition is for finite dimensional matrices and needs to
modified when $\Abf$ is a closed linear operator on a sub-space of a
Hilbert-space. In this case, we extend these definitions
by replacing ${\bf u}$ with a sequence of functions, $\{ {\bf u}_n \}$, in the
domain of  $\Abf$, i.e.~require that Definition 1 hold on
the closure of the  domain of $\Abf$. Thus, condition (iii) becomes
$||({\bf A}-\lambda{\bf I})^{-1} || \geq 1/\epsilon$.

In [8,10,11], two methods for computing the $\epsilon$-pseudospectrum 
are given.
First, using (i), we define $\epsilon_b (z)$ to be the 
smallest singular value of ${\bf A}-z{\bf I}$.
The subscript $b$ on $\epsilon$
is used because $\epsilon_b (z)$
is the boundary of the $\epsilon$-pseudospectrum: $z \in
\Lambda_{\epsilon} ({\bf A})  \Leftrightarrow  \epsilon_b (z) \leq \epsilon$.

$\epsilon_b (z)^2$ can also be determined by computing the 
smallest eigenvalue of $({\bf A}-z{\bf I})^{*}({\bf A}-z{\bf I})$,
which is equivalent to minimizing
$\parallel ({\bf A}-z{\bf I}){\bf u} \parallel^2$ over ${\bf u}$
with $||{\bf u}||^2 = 1$. However, this alternative is not as well 
conditioned as  the singular value decomposition of ${\bf A}-z{\bf I}$
due to rounding error computing $({\bf A}-z{\bf I})^{*}({\bf A}-z{\bf I})$.

Second, for a fixed value of $\epsilon$,
Trefethen and Reddy et al.~approximately compute the 
$\epsilon$-pseudospectrum by
generating random perturbing matrices, $\bf E$, of unit norm, and then
calculating the eigenvalues, $\{z({\bf E}, \epsilon )\}$, of
${\bf A}+ \epsilon {\bf E}$. The eigenvalues, $\{z({\bf E}, \epsilon )\}$, 
are plotted in the complex $z$-plane for the ensemble of matrices, 
${\bf E}$. As the number of perturbing matrices increases, the scatterplot of
the calculated eigenvalues densely fills the $\epsilon$-pseudospectrum
(provided that the ensemble is representative of all possible perturbations).

We can modify this algorithm by considering only rank 1 perturbing matrices,
${\bf E}$. This approach has several advantages. First, the space of all rank
one perturbations has dimension $2(n-1)$, while the space of all matrices
with unit norm has dimension $n^2 -1$. Thus, for a given number, $N$, of
perturbing matrices, the distance, as measured in trace norm,
between an arbitrary element, ${\bf E}^{\prime}$, and the closest matrix
in the test ensemble is $O(N^{-1/d})$, where $d$ equals $2(n-1)$ and   
$n^2 -1$ respectively. Thus the typical distance to a matrix in the 
test ensemble is significantly smaller in the rank one case.
Using only rank 1 perturbations of the form,
$({\bf Av} - z {\bf v}) {\bf v}^{*}$, further reduces the space of possible
perturbations.

Second, for general perturbing matrices, $\bf E$, a singular value
decomposition needs to be computed to determine $\parallel {\bf E} \parallel$.
This step is unnecessary for rank one perturbations. 

Third, if the matrix
has a special structure, an arbitrary perturbation will
destroy these characteristics. 
In contrast, rank one perturbations only modify the special structure.
In certain cases, it may be possible to use this special structure
in specialized algorithms to compute the eigenvalue decomposition.

Another important stability question  is  
``{\it How close is a matrix, $\Abf$, to being unstable?}\ ''
In Ref.~13, Van Loan addressed this question. We now reformulate his
results in the language of $\epsilon$-psuedospectra.

\begin{definition}
{\it Let $\epsilon \ge 0$ and $\Abf \in C^{m \times m}$, 
$\Abf$ is $\epsilon$-asymptotically stable if and only if there exists 
no  matrix, $\Ebf \in C^{m \times m}$ such that
$||\Ebf|| < \epsilon$ and $\Abf + \Ebf$ has an eigenvalue, $\lambda$,
with $\Re[\lambda] \ge 0$.}
\end{definition}

\begin{theorem}
{\it 
$\Abf$ is $\overline{\epsilon}$-asymptotically stable if and only if
$\inf_{y \in R} \eps_b(iy) \ge \overline{\epsilon}$.
}
\end{theorem}

Thus the minimum of $\eps_b(z)$ on the imaginary axis defines the 
$\eps$-instability threshold.

\section{Generalized singular value decomposition}

Before examining generalized pseudospectra, we  review
the generalized singular value decomposition [1,7,12].

\begin{definition}
{\rm (Van Loan)}. {\it
The $\Bbf$-{\it singular values} of a matrix $\Abf$ are elements of the set
$\mu (\Abf,\Bbf)$ defined by
$$
\mu (\Abf,\Bbf) \equiv \{ \mu | \mu \ge 0, \ {\rm det} \ (\Abf^{*}
\Abf - \mu^2 \Bbf^{*} \Bbf) = 0 \} \ .
$$
where $\Abf \in C^{m_a \times n}$, $\Bbf \in C^{m_b \times n}$ and 
$m_a \ge n$.} 
\end{definition}

Van Loan's original work decomposed real matrices with orthogonal
matrices. We state the analogous result for complex matrices and 
unitary transformations. 

\begin{theorem}
{\rm (Van Loan) (The $\Bbf$-singular value decomposition (BSV).)}
{\it Suppose
$\Abf \in C^{m_a \times n}$, $\Bbf \in C^{m_b \times n}$ and $m_a \ge n$. There
exist unitary matrices, $\Ubf (m_a \times m_a )$ and $\Vbf (m_b \times m_b )$, 
and a nonsingular $n \times n$ matrix, $\Xbf$, such that
$$
\Ubf^{*} \Abf\Xbf = \Dbf_{\Abf} = 
{\rm diag} \ ( \alpha_1 , \ldots , \alpha_n ) \ ,
\ \ \ \ \alpha_i \ge 0 \ ,
$$
$$
\Vbf^{*} \Bbf\Xbf = \Dbf_{\Bbf} = {\rm diag} \ ( \beta_1 , \ldots , \beta_n ) \ ,
\ \ \ \ \beta_i \ge 0 \ ,
$$
where $q = \min \{ m_b , n \}$, $r =$ rank$(\Bbf)$ and $\beta_1
\ge \ldots \ge \beta_r > \beta_{r+1} = \ldots =
\beta_q = 0$. If $\alpha_j = 0$ for any $j$, $r+1 \le
j \le n$, then $\mu (\Abf,\Bbf) = \{ \mu | \mu
\ge 0 \}$. Otherwise, $\mu (\Abf,\Bbf) = \{ \alpha_i / \beta_i |
i = 1, \ldots , r \}$.
}
\end{theorem}

In [7], Paige and Saunders generalized this definition to relax the
requirement that $m_a \ge n$. 
Van Loan also gave a second generalization of the singular value decomposition,
by using two different norms in the variational formulation.

\begin{definition}
{\rm (Van Loan).}
{\it Let $\Pbf \in C^{n \times n}$ be positive definite.
A matrix
$\Qbf \in C^{n \times n}$ is $\Pbf$ unitary if 
$\Qbf^{*} \Pbf\Qbf = \Ibf_n$.}
\end{definition}

\begin{definition}
{\rm (Van Loan).} {\it
Let $\Sbf$ and $\Tbf$ be positive definite matrices of orders $m$
and $n$, respectively, with $m \ge n$. The {\rm(}$\Sbf,\Tbf${\rm)}-singular
values of $\Abf \in C^{m \times n}$ are elements of 
the set $\mu (\Abf,{\Sbf,\Tbf})$ defined by }
\end{definition}
$$
\mu(\Abf,{\Sbf,\Tbf}) = \{ \mu | \mu \ge 0, \ \mu^2 \ {\rm is \ a \
stationary \ value \ of} \ 
{\xbf^*\Abf^* \Sbf\Abf \xbf \over \xbf^*\Tbf\xbf } \} \ .
$$

\begin{theorem}
{\rm (Van Loan) (The ($\Sbf,\Tbf$)-singular value decomposition}{\rm)}.
{\it Let $\Abf,\Sbf$ and $\Tbf$ be
in $C^{m \times n}$, $C^{m \times m}$, $C^{n \times n}$, respectively, 
with $\Sbf$ and $\Tbf$
positive definite $(m \ge n)$. There exists an $\Sbf$ unitary
$\Ubf \in C^{m \times m}$ and a $\Tbf$ unitary $\Vbf \in C^{n \times n}$ such
that}
\end{theorem}
$$\Ubf^{-1} \Abf\Vbf = \Dbf = {\rm diag} \ ( \mu_1 , \ldots , \mu_n )
\ \ .
$$
We denote the ($\Sbf,\Tbf$)-singular values of $\Abf$ by
$\mu(\Abf,\ \Sbf,\ \Tbf)$.
The proof is based on the singular value decomposition of 
$\Lbf^{*} \Abf\Kbf^{-1^*}$, where $\Kbf$ and $\Lbf$ are the
Cholesky factorizations of $\Tbf$ and $\Sbf$: $\Tbf = \Kbf\Kbf^{*}$ 
and $\Sbf = \Lbf\Lbf^{*}$. 
We relate the two generalizations of the singular value decomposition 
by the following corollary.

\begin{cor} 
{\rm (The ($\Sbf,\Tbf$)-singular value decompositions
corresponding to the $\Bbf$-singular value decomposition).} 
{\it Let $\Abf$ and $\Bbf$ be in $C^{m \times m}$,  with $\Bbf$ 
positive definite, self-adjoint. Then}
\end{cor} 

{\it 
a)\ $\mu (\Abf,\Bbf) = \mu (\Abf ,\ \Ibf_m,\ \Bbf^{*}\Bbf)$


b)\ $\mu (\Abf,\Bbf) = \mu(\Abf ,\ {\Bbf^{-1 }\Bbf^{-1^*}, \ \Ibf_m})$
}\ .

The corollary follows from determinant identities. 

\section{Generalized $\epsilon$-pseudospectrum}

We now consider the generalized eigenvalue problem
${\bf A}{\bf e} = \lambda {\bf M}{\bf e}$, where
$\bf M$ is self-adjoint and positive definite. We could transform the problem
into a standard eigenvalue problem: ${\bf A}_{\bf F}
{\bf e}^{\prime} =  \lambda{\bf e}^{\prime}$
where ${\bf F}^{*} {\bf F} = {\bf M}$,  ${\bf e}' = {\bf F}{\bf e}$, and
${\bf A_F} \equiv{\bf F}^{-1^* } {\bf A}{\bf F}^{-1}$,
and then examine the $\epsilon$-pseudospectrum
of the standard linear problem. 

In Appendix A of [8], Reddy et al.~consider the $\epsilon$-pseudospectrum 
of an  operator, $\Abh$, where both the domain and the range of $\Abh$
have a metric $\Mbf= \Fbf^*\Fbf$. In this case, the $\epsilon$-pseudospectrum 
is defined as the set of $z$ values such that  there exists a vector
$\ubf$ with ${\bf u}^{*} {\Abh}^{*}
{\bf M} {\Abh}{\bf u} \leq \eps^2$ and  ${\bf u}^{*} {\bf M}{\bf u} = 1$.
Reddy et al.~show that the $\epsilon$-pseudospectrum for the
$ \Fbf^*\Fbf$ norm is equivalent to the standard $\epsilon$-pseudospectrum 
for $ \Fbf\Abh\Fbf^{-1}$. 

We modify Reddy et al.'s approach in two way. First, we consider
the generalized eigenvalue problem with
${\bf A}{\bf e} =  \lambda {\bf M}{\bf e}$ with 
$\Abf \equiv \Fbf^*\Fbf \Abh$  and  $\Mbf= \Fbf^*\Fbf$. This type
of  generalized eigenvalue problem  
tends to occur when differential equations are discretized using a 
variational formulation. For this formulation, the $\epsilon$-pseudospectrum 
of the generalized eigenvalue problem is equivalent to the standard 
definitions of $\epsilon$-pseudospectrum for $ \Fbf^{-1^*}\Abf\Fbf^{-1}$.   

Second, we prefer not to transform the generalized $\epsilon$-pseudospectrum 
problem into an equivalent standard problem. 
The transformation ${\bf A} \rightarrow  \Fbf^{-1^*}\Abf\Fbf^{-1}$ 
disguises the effect of the matrix norm and 
can result in the loss of accuracy due to roundoff error.
For similar problems such as the generalized singular value decomposition and
the generalized eigenvalue decomposition, the best numerical algoriithms do not
transform the generalized problem, but solve it directly. 
Therefore, we state all of the equivalent definitions of
$\epsilon$-pseudospectra for the generalized case.
We do this for our formulation, with $\Abf \equiv \Fbf^*\Fbf \Abh$.
The Reddy et al.~formulation is given by transforming the results
in definition 3.1.

We restrict our consideration to the finite dimensional case.   
We denote the spectrum of the generalized eigenvalue problem,
${\bf A}{\bf e} = \lambda {\bf M}{\bf e}$, by
$\Lambda({\bf A}, {\bf M})$ and 
the resolvent set by $\rho({\bf A},{\bf M})$.

\begin{definition}
 {\it
Let $\Mbf$ be a positive self-adjoint matrix and let
$\epsilon \geq 0$ be given. A complex number $z$ is in the
$\epsilon$-pseudospectrum of $( {\bf A,M})$,
which we denote by $\Lambda_{\epsilon} ( {\bf A,M})$,
if any of the following equivalent conditions is satisfied:

(0) z is in the $\epsilon$-pseudospectrum of $ \Fbf^{-1^*} {\bf AF}^{-1}$,
where $\Fbf^{*} \Fbf = {\bf M}$.

(i) the smallest generalized $({\bf M}^{-1}\ , \ \Mbf)$  
singular value of ${\bf A}-z{\bf M}$
is less than or equal to $\epsilon$,
i.e. $\epsilon \ge min\{\mu({\bf A}-z{\bf M},\Mbf^{-1},\Mbf) \}$.

(ii) $\exists {\bf u}\in C^n$ such that
${\bf u}^{*} {\bf {\bf Mu}} = 1$ and ${\bf u}^{*} ({\bf A}-z{\bf M})^{*}
{\bf M}^{-1} ({\bf A}-z{\bf M}) {\bf u} \leq \epsilon^2$,

(iii) $z \in \Lambda ({\bf A,M})$ or
$z \in \rho ({\bf A,M})$ and $\exists {\bf w} \in C^n$
such that ${\bf w}^{*} {\bf M}^{-1}{\bf w} = 1$
and ${\bf w}^{*} ({\bf A}-z{\bf M})^{-1^*}$
${\bf M}({\bf A}-z{\bf M})^{-1} {\bf w} \geq 1/
\epsilon^2$,

(iv) $z$ is in the generalized spectrum of
${\bf A}+ \epsilon \Fbf^{*} {\bf E}\Fbf$:
$({\bf A} + \epsilon
\Fbf^{*} {\bf E}\Fbf) {\bf u} = z{\bf M}{\bf u} $,
where $\Fbf^{*} \Fbf = {\bf M}$ and the matrix ${\bf E}$ satisfies
$\parallel {\bf E} \parallel \leq 1$,

(iv')
$\exists$ a $n \times n$ matrix, ${\bf H}$, such that
$z$ is in the generalized spectrum of ${\bf A}+ \epsilon {\bf H}:
({\bf A}+ \epsilon {\bf H}) {\bf u}
=z{\bf Mu}$, where the matrix ${\bf H}$ satisfies
$$
\max_{{\bf u} \in C^n} {{\bf u}^{*} {\bf H}^{*}
{\bf M}^{-1} {\bf H}{\bf u} \over {\bf u}^{*} {\bf M}{\bf u}} \leq 1
\ \ , $$

(v) $\exists {\bf u}_1$ and ${\bf u}_2 \in C^n$ such that $z$ is in the 
generalized spectrum of ${\bf A} + \epsilon {\bf w}_2
{\bf w}_1^{*}$ w.r.t. {\bf M}, where ${\bf w}_1 = {\bf Mu}_1$,
${\bf w}_2 = {\bf Mu}_2 ,$ ${\bf u}_1^{*} {\bf Mu}_1 \leq 1$
and ${\bf u}_2^{*} {\bf Mu}_2 \leq 1$.

(vi) $\exists {\bf u}$  in $C^n$ such that $z$ is in the generalized
spectrum of ${\bf A} -  ( {\bf Au}-z{\bf Mu} )
{\bf w}^{*}$ w.r.t. {\bf M}, where ${\bf w} = {\bf Mu}$,
${\bf u}^{*} {\bf Mu} \leq 1$
and ${\bf u}^{*} ({\bf A}-z{\bf M})^{*}
{\bf M}^{-1} ({\bf A}-z{\bf M}) {\bf u} \leq \epsilon^2$.
}
\end{definition}

The equivalence may be proved directly or by simply transforming each
of the properties from definition 1.1 as applied to ${\bf A}_\Fbf$.
Properties (i) and (v) are used in practice. Property (ii) corresponds
to part b) of the corollary. \hfill $\Box$

We now present a different  generalization of
$\epsilon$-pseudospectra for the generalized eigenvalue problem,
$\Abf \ubf =\lambda \Mbf \ubf$. Our $\Mbf$-weighted $\epsilon$-pseudospectrum
has the advantage that definitions (i)-(vi) are simpler than in
Def.~3.1. However, the $\Mbf$-weighted $\epsilon$-pseudospectrum
is not related to the standard $\epsilon$-pseudospectrum
of Def.~1.1 through a change of variables, and therefore Def.~3.2 has no
analog of (0) in Def.~3.1. 

\begin{definition}
{\rm ($\Mbf$-weighted $\epsilon$-pseudospectrum).} {\it
Let $\Mbf$ be a positive self-adjoint matrix and 
let $\epsilon \geq 0$ be given. Define $\epsbr \equiv \eps 
\left( {||M|| \over ||M^{-1}||} \right)^{1\over 2}$.
A complex number $z$ is in the
$\Mbf$-weighted $\epsilon$-pseudospectrum of ${\bf A}$,
which we denote by $\Lambda_{\epsilon}( {\bf A} | \Mbf )$,
if one of the following equivalent conditions is satisfied:

(i) the smallest 
singular value of ${\bf A}-z{\bf M}$ is
less than or equal to $\epsbr$.

(ii) $\exists {\bf u}\in C^n$ such that
$||{\bf u}||^2 = 1$ and 
$ ||({\bf A}-z{\bf M}) {\bf u}||^2 \leq \epsbr^2$,

(iii) $z \in \Lambda ({\bf A})$ or
$z \in \rho ({\bf A})$ and $\exists {\bf u} \in C^n$
such that $||{\bf u}||^2 = 1$
and ${\bf u}^{*} ({\bf A}-z{\bf M})^{-1 *}$
$({\bf A}-z{\bf M})^{-1} {\bf u} \geq 1/
\epsbr^2$,

(iv) $z$ is in the generalized spectrum of
${\bf A}+ \epsbr  {\bf E}$ w.r.t. {\bf M}:
$({\bf A} + \epsbr{\bf E}) {\bf u} = z{\bf Mu} $,
 where the matrix ${\bf E}$ satisfies
$\parallel {\bf E} \parallel \leq 1$.

(v) $z$ is in generalized spectrum of 
${\bf A} + \epsbr {\bf v}_2 {\bf v}_1^{*}$
w.r.t. ${\bf M}$, where
$\parallel {\bf v}_1 \parallel \leq1$ and
$\parallel {\bf v}_2 \parallel \leq1$.

(vi) $z$ is in the generalized spectrum of 
${\bf A} - ({\bf Av} - z {\bf Mv}) {\bf v}^{*}$
w.r.t. ${\bf M}$, where
$\parallel {\bf v} \parallel \leq1$ and
$\parallel {\bf Av} - z {\bf Mv} \parallel \leq \epsbr$.
}
\end{definition}

The normalization, $\epsbr \equiv \eps 
\left( {||\Mbf|| / ||\Mbf^{-1}||} \right)^{1\over 2}$,
allows Def.~3.2 to reduce to Def.~1.1 when $\Mbf$ is a multiple of the 
identity matrix. When $||\Mbf^{-1}||$ in infinite, we can replace this
normalization with the normalization: $\epsbr \equiv \eps ||\Mbf||$.
When $\Mbf$ is an unbounded operator with finite $||\Mbf^{-1}||$,
we can replace this definition with the nomalization: 
$\epsbr \equiv \eps / ||\Mbf^{-1}||$.

Since definition 3.1 is a transformed version of Def.~1.1, we believe
that Def.~3.1 is preferable to the simpler, but  coordinate dependent 
Def.~3.2. In particular,
definition 3.1 is useful in the analysis of differential
operators which have variational formulations.

In definition 3.1, $\Mbf$ is the metric of the domain of $\Abf$
and  $\Mbf^{-1}$ is the metric of the range of $\Abf$. 
In the corresponding Reddy et al.~formulation, $\Mbf$ is the metric 
of both the domain and the range of $\Abh$. We give a more general formulation
of $\epsilon$-pseudospectra which incorporates both previous cases.
We say that $z$ is in the
$\epsilon$-pseudospectrum of $( {\bf A,M})$ with respect to the operator norms
$\parallel \cdot \parallel_1$ and $\parallel \cdot \parallel_2$ if and only if
there exist matrices $\Ebf_1$ and $\Ebf_2$ such that 
$\parallel \Ebf_1 \parallel_1^2 + \parallel \Ebf_2 \parallel_2^2 \leq \eps^2$
and $z$ is a  generalized eigenvalue of $( {\Abf+\Ebf_1 ,\Mbf +\Ebf_2})$.
This extended
definition corresponds to part (iv) of Defs.~1.1 \& 3.1. 
Definition 3.1 
is equivalent  
to restricting to $\Ebf_2 \equiv 0$ and using the matrix norm: 
$\parallel \Ebf \parallel_{\Mbf,\Mbf^{-1}}^2 \equiv
\max_{{\bf u} \in C^n} {{\bf u}^{*} {\bf E}^{*}
{\bf M}^{-1} {\bf E}{\bf u} \over {\bf u}^{*} {\bf M}{\bf u}}$.

Restricting to $\Ebf_2 \equiv 0$ is natural when $\Mbf$ is known to
higher precision than $\Abf$. When $\parallel \cdot \parallel_1 \equiv
\parallel \cdot \parallel_2$, then allowing for both
$\Ebf_1$ and $\Ebf_2$  simply transforms the definitions:
$\Ebf_1= {\Ebf \over 1+|z|^2}$ and $\Ebf_2 = {
-\overline{z} \  \Ebf \over 1+|z|^2}$, where $\Ebf$ is the optimal 
perturbation. The resulting perturbation has norm:
${\parallel \Ebf_1 \parallel^2 + \parallel \Ebf_2 \parallel^2 }
= {1 \over 1+|z|^2} \parallel \Ebf \parallel^2$. Thus by allowing an
$\Ebf_2$ perturbation, we  replace the critical value of $\eps$, $\eps_b(z)$,
with ${1 \over \sqrt{1+|z|^2}}\eps_b(z)$.
J.~Demmel points out that a modification of the proof of Lemma 5 of [3] 
shows the  equivalence of (i) and (iv) in Def.~3.2.

\section{Transient growth for the generalized system}

For completeness, we now describe  algorithms for computing the numerical
range and maximum transient growth for the generalized eigenvalue
problem.

\begin{definition}
{\it
Let ${\bf A}$ be a linear operator and ${\bf M}$ be a self-adjoint,
positive definite operator. 
The numerical range of $\bf A$ with respect to $\bf M$ is defined by }
$$
R({\bf A,M}) \equiv \{ z|\ \exists {\bf u} \ {\rm with } \
{\bf u}^{*} {\bf A} {\bf u} = z \ {\rm and} \
{\bf u}^{*} {\bf Mu} = 1 \}
\ \ .$$
\end{definition}

\noindent
The numerical range is convex, and its boundary
can be computed by maximizing $Re (e^{- i \theta }
\lambda_{\theta} )$ for all $\theta$, 
where $\lambda_{\theta}$ is a generalized eigenvalue.
Thus for each value of $\theta \in [- \pi , \pi ]$, the largest
eigenvalue, $\lambda_{\theta}$, of the self-adjoint eigenvalue problem is computed,
$$
(e^{- i \theta } {\bf A} + e^{i \theta }
{\bf A}^{*} ){\bf e}_{\theta} = 2 \lambda_{\theta}
{\bf M}{\bf e}_{\theta}
\ \ .$$
We note that the $\epsilon$-psuedospectrum is contained within $\epsilon$
of the numerical range [7,9,10].  
\vspace{.3in}

Finally, we consider transient growth problems for evolutionary systems
of  partial differential equations. In [8], Reddy et al.~studied 
the transient growth of solutions of the Orr-Sommerfeld equation.
The corresponding problem in magneto-hydrodynamics again requires a 
generalized eigenvalue problem. We consider 
the initial value problem:
$$
{\bf M} {\partial {\bf u} \over \partial t} = {\bf Au}
\ \ .$$

We wish to maximize the energy at time $t$, $E_{\Mbf} (t) \equiv
{\bf u}(t)^{*} {\bf Mu}(t)$, subject
to $E_{\Mbf} (t=0 ) = 1$. We compute the eigenvalues and eigenvectors,
$\{ ( \lambda_k ,{\bf e}_k ) \}$ of ${\bf Ae}_k = \lambda_k {\bf Me}_k$. 
We represent
${\bf u}(t) = \sum_k a_k {\bf e}_k e^{\lambda_k t} $. We define
$Q_{k, \ell} (t) =
{\bf e}_k^{*} {\bf Me}_{\ell} \exp^{( \bar{\lambda}_k + \lambda_{\ell} )t}
$.
The energy at time, $E_{\Mbf} (t)$, is
${\bf a}^{*} {\bf Q}(t){\bf a}$,
where $\bf a$ is the $n$ vector of coefficients.
Thus the maximum transient growth is computed by maximizing
$E_{\Mbf}(t)- \lambda (E_{\Mbf} (t=0) -1)$ with respect to ${\bf a}$.
The resulting self-adjoint eigenvalue problem is
$$
{\bf Q}(t) {\bf a} (t) = \lambda (t) {\bf Q}(t = 0) {\bf a}(t)
\ \ ,$$
where ${\bf a}(t)$ are the coefficients which maximize the transient growth,
and $\lambda$ is the maximium energy growth. 
Transient growth depends on  the norm which measures the energy.
Replacing  $E_{\Mbf} (t)$ with $E_{\bf N} (t) \equiv
{\bf u}(t)^{*} {\bf Nu}(t)$, where ${\bf N}$ is an arbitrary positive definite
matrix, can   greatly alter the magnitude of the transient growth. 
This maximum transient growth problem can  be recasted
as a generalized singular value problem for $\mu (\Abf,\Bbf)$ where
${\bf F}^* {\bf F} = {\bf M}$,
$\Bbf_{\cdot,\ell} \equiv \Fbf \ebf_{\ell}$, \
$\Abf_{\cdot,\ell} \equiv e^{\lambda_{\ell} t}\Fbf \ebf_{\ell}$,
and $\mu^2(t) = \lambda(t)$.
This generalized singular value decomposition formulation 
reduces the roundoff error in the computation.

In conclusion, $\epsilon$-pseudospectra, the numerical range, and the maximum
transient growth rate have been useful in analyzing 
certain problems in fluid dynamics. For evolutionary systems, the
$\epsilon$-pseudospectrum describes the norm of the Green's function to 
fixed frequency forcing.
To treat magneto-hydrodynamics, the corresponding definitions and algorithms
for the generalized eigenvalue problem are required. Parts (v) and (vi) of
definitions 1.1 and 3.1 simplify the calculation of $\epsilon$-pseudospectra.

\medskip

\end{document}